\documentclass[a4paper, 12pt]{article}
\usepackage{amsfonts}
\usepackage{fancyhdr}
\usepackage{graphicx}
\usepackage{epsfig}
\usepackage{amssymb}
 \usepackage{mathrsfs,amsfonts,amsmath}
 \usepackage{color}

\makeatletter

\@addtoreset{equation}{section} \makeatother
\newtheorem{Theorem}{Theorem}[section]
\newtheorem{Remark}[Theorem]{Remark}
\newtheorem{Lemma}[Theorem]{Lemma}

\begin{document}
\title{\textbf{Strong convergence rates of modified truncated EM methods for neutral stochastic differential delay equations
\footnote{Supported by Natural Science Foundation of China (NSFC
11601025).}}}
\author{ Guangqiang Lan\footnote{Corresponding author: Email:
langq@mail.buct.edu.cn.}\quad and\quad
 Qiushi Wang
\\ \small School of Science, Beijing University of Chemical Technology, Beijing 100029, China}

\date{}

\maketitle

\begin{abstract}
The aim of this paper is to investigate strong convergence of modified truncated Euler-Maruyama method for neutral stochastic differential delay equations introduced in Lan (2018). Strong convergence rates of the given numerical scheme to the exact solutions at fixed time $T$ are obtained under local Lipschitz and Khasminskii-type conditions. Moreover, convergence rates over a time interval $[0,T]$ are also obtained under additional polynomial growth condition on $g$ without the weak monotonicity condition (which is usually the standard assumption to obtain the convergence rate). Two examples are presented to interpret our conclusions.
\end{abstract}

\noindent\textbf{MSC 2010:} 60H10, 65C30, 65L20.

\noindent\textbf{Key words:} neutral stochastic differential delay equations; local Lipschitz condition; modified truncated Euler-Maruyama method; strong convergence rate.

\section{Introduction}

\noindent

Let $(\Omega,\mathscr{F},\{\mathscr{F}_t\}_{t\geq 0},P)$ be a complete
filtered probability space satisfying usual conditions (i.e. $\{\mathscr{F}_t\}_{t\geq 0}$ is right continuous and $\mathscr{F}_0$ contains all $P$-null sets). Consider the following neutral stochastic differential delay equations (short for NSDDEs)

\begin{equation}\label{nsdde}d(x(t)-D(x(t-\tau))
=f(x(t),x(t-\tau))dt+g(x(t),x(t-\tau))dB_t
\end{equation}
on $t\ge0$ with initial value $\{x(\theta):-\tau\le\theta\le0\}=\xi\in C_{\mathscr{F}_0}^b([-\tau,0],\mathbb{R}^d)$ (the family of all $\mathscr{F}_0$ measurable bounded $C([-\tau,0],\mathbb{R}^d)$-valued random variables), for any $\xi\in C_{\mathscr{F}_0}^b([-\tau,0],\mathbb{R}^d), p\ge1$ define  $||\xi||_p=(\mathbb{E}(\sup_{-\tau\le\theta\le0}|\xi(\theta)|^p))^{\frac{1}{p}}$, $\{B_t,t\geq0\}$ is an $n$-dimensional standard $\mathscr{F}_t$-Brownian motion, $D:x\in\mathbb{R}^d\mapsto D(x)\in\mathbb{R}^d$, $f:(x,y)\in\mathbb{R}^d\times\mathbb{R}^d\mapsto f(x,y)\in\mathbb{R}^d$ and
$g:(x,y)\in\mathbb{R}^d\times\mathbb{R}^d\mapsto g(x,y)\in\mathbb{R}^d\otimes\mathbb{R}^n$ are measurable functions.

Recently, such neutral stochastic differential delay equations (short for NSDDEs) have been found more and more applications in many fields such as control theory, electrodynamics, biomathematics and so on. However, most NSDDEs can not be solved explicitly except some special ones. Thus numerical methods for NSDDEs (\ref{nsdde}) have been playing more and more important roles.

The convergence of the numerical methods for NSDDEs (\ref{nsdde}) have been discussed intensively by many researchers, for example, Gan et. al \cite{GSZ} investigated mean square convergence of stochastic $\theta$ method under global Lipschitz condition, \cite{ZG} studied Mean square convergence of one-step methods under the same assumptions, \cite{M} considered convergence in probability of the backward Euler approximate solution for a class of stochastic differential equations with constant delay, Zhang et.al considered strong convergence of the partially truncated Euler-Maruyama method for a class of stochastic differential delay equations in \cite{ZSL}, there are also many other literatures concerning with this topic, see e.g. \cite{EXT,Kloeden,mao2,MY,ZSL1,ZW}.

Recently, in \cite{mao}, Mao developed a new explicit numerical simulation method, called truncated EM method. Strong convergence theory were established there under local Lipschitz condition plus the Khasminskii-type condition. And then he obtained sufficient conditions for the strong convergence rate of it in \cite{mao1}. Motivated by these two works, Lan and Xia introduced in \cite{LX} modified truncated Euler-Maruyama (MTEM) method and obtained the strong convergence rate under given conditions. Then in \cite{L}, the author generalized the MTEM method from SDEs cases to the NSDDEs cases and obtain asymptotic exponential stability of it under given conditions. However, the strong convergence of the MTEM method is still not known, which is the main topic of this paper.

Although strong convergence of the given numerical methods are considered in many papers such as \cite{BT,EXT,GSZ,mao,YM,ZG} and so on, the convergence rates are not known. Some other papers considered strong convergence rates of the given numerical methods under weak monotonicity condition. For example, in Guo et al \cite{GMY}, Assumption 5.1 is necessary to obtain the convergence rate of the truncated EM method at fixed time $T$, in Tan and Yuan \cite{TY}, A8 is needed to obtain the convergence rate of the truncated EM method over a time interval $[0,T]$, for the weak monotonicity, one can also see \cite{mao1,LX,TY1} and so on. In this paper, we will consider the strong convergence rates of MTEM methods to exact solutions both at fixed time $T$ and over a time interval $[0,T]$ without such weak monotonicity conditions.

The organization of the paper is as the following. In Section 2, the MTEM method for NSDDE is introduced, and main results are presented. In Section 3, some useful lemmas are presented to prove the convergence theorems. In Section 4, convergence rates at fixed time $T$ are obtained. The convergence rates over the time interval $[0,T]$ will be proved with additional polynomial growth condition on $g$ in Section 5. Then in Section 6, two examples are presented to interpret the Theorems. We will conclude our paper in Section 7.

\section{The settings and main results}

Assume that both the coefficients $f$ and $g$ in (\ref{nsdde}) are locally Lipschitz continuous, that is, for each
$R$ there is $L_R>0$ (depending on $R$) such that
\begin{equation}\label{local}|f(x,y)-f(\bar{x},\bar{y})|
\vee|g(x,y)-g(\bar{x},\bar{y})|\le L_{R}(|x-\bar{x}|+|y-\bar{y}|)\end{equation} for all
$|x|\vee|y|\vee|\bar{x}|\vee|\bar{y}|\le R>0$. Here the norm of a matrix $A$ is denoted by $|A|=\sqrt{\textrm{trace}(A^\textrm{T}A)}$.

Assume also that there is a positive constant $u\in(0,1)$ such that

\begin{equation}\label{Lip}|D(x)-D(y)|\le u|x-y|,\ \forall x,y\in\mathbb{R}^d\end{equation}

It is well known that there is a unique strong solution (might explode at finite time) to equation (\ref{nsdde}) under conditions (\ref{local}) and (\ref{Lip}), see e.g. \cite{MY}.

As interpreted in \cite{LX}, we can always choose $\Delta^*>0$ small enough and a strictly positive decreasing function $h:(0,\Delta^*]\to(0,\infty)$ such that
\begin{equation}\label{tiaoj}\lim_{\Delta\to0}h(\Delta)=\infty\ \textrm{and}\ \lim_{\Delta\to0}L_{h(\Delta)}^4\Delta=0.\end{equation}

For any $\Delta\in(0,\Delta^*),$ we define the the modified truncated function of $f$ as the following:
\begin{equation}\label{dy}f_\Delta(x,y)=\left\{\begin{array}{ll} f(x,y),\qquad\qquad\qquad\ |x|\vee|y|\le h(\Delta),\\
\frac{|x|\vee|y|}{h(\Delta)} f\left(\frac{h(\Delta)}{|x|\vee|y|}(x,y)\right),|x|\vee|y|>h(\Delta).
\end{array}
\right.\end{equation}
$g_\Delta$ is defined in the same way as $f_\Delta$. Here $f(a(x,y))\equiv f(ax,ay)$ for any $a\in \mathbb{R}, x,y\in\mathbb{R}^d.$

It is obvious that the functions $f_\Delta$ and $g_\Delta$ defined above are different from the truncated functions defined in \cite{GMY}.

We have defined the discrete MTEM method in \cite{L}. However, we recall it here for readers' convenience.

Let $\Delta$ be a stepsize such that $\tau=m\Delta$ for some positive integer $m$. Then by using $f_\Delta$ and $g_\Delta$, we can define the MTEM method $X_k$ of (\ref{nsdde}) as the following:
\begin{equation}\label{num}\aligned X_{k+1}-D(X_{k+1-m})&=X_k-D(X_{k-m})+f_\Delta(X_k,X_{k-m})\Delta\\&\quad
+g_\Delta(X_k,X_{k-m})\Delta B_k,\quad k=0,1,2,\cdots,\\&
X_k=\xi(k\Delta),\quad k=-m,-m+1,\cdots,0.\endaligned\end{equation}
Here $\Delta B_k=B((k+1)\Delta)-B(k\Delta)$ is the increment of the $n$-dimensional standard Brownian motion.

The two versions of the continuous-time MTEM solutions are defined as the following:

\begin{equation}\label{num1}\bar{x}_\Delta(t)=\sum_{k=-m}^\infty X_k1_{[k\Delta,(k+1)\Delta)}(t),\quad t\ge-\tau,\end{equation}
and $x_\Delta(t)=\xi(t),$ $t\in[-\tau,0],$
\begin{equation}\label{num2}\aligned x_\Delta(t)&=D(\bar{x}_\Delta(t-\tau))+\xi(0)-D(\xi(-\tau))+\int_0^tf_\Delta(\bar{x}_\Delta(s),\bar{x}_\Delta(s-\tau))ds\\&
\quad+\int_0^tg_\Delta(\bar{x}_\Delta(s),\bar{x}_\Delta(s-\tau))dB(s),\quad t\ge0.\endaligned\end{equation}

Obviously, $x_\Delta(k\Delta)=\bar{x}_\Delta(k\Delta)=X_k$ for all $k\ge0.$

To study the strong convergence of continuous version of MTEM (\ref{num2}), let us first consider the following condition:

Assume that there exist positive constants $K$ and $p>2$ such that
\begin{equation}\label{dandiao}
2\langle
x-aD(\frac{1}{a}y),f(x,y)\rangle+(p-1)|g(x,y)|^2\le K(1+|x|^2+|y|^2)
\end{equation}
holds for all $x,y\in\mathbb{R}^d, a\in(0,1].$

Notice that when $a=1,$ (\ref{dandiao}) reduces to the well known Khasminskii condition
\begin{equation}\label{dandiao1}
2\langle
x-D(y),f(x,y)\rangle+(p-1)|g(x,y)|^2\le K(1+|x|^2+|y|^2).
\end{equation}

Suppose for fixed $\Delta$ ($\tau=m\Delta$) the initial value $\xi$ satisfies
\begin{equation}\label{chuzhi}\mathbb{E}\sup_{-m\le k\le -1}\sup_{k\Delta\le s\le (k+1)\Delta}|\xi(s)-\xi(k\Delta)|^q\le \hat{K}\Delta^\frac{q}{2}\end{equation}
for $2<q<p$.

Now we are ready to state our first result on the strong convergence rate for MTEM method at fixed time $T.$

\begin{Theorem}\label{conv}
Assume that (\ref{local}), (\ref{Lip}), (\ref{dandiao}) and (\ref{chuzhi}) hold for some $2<q<p$. If there exist $0<\Delta_0\ (\le\Delta^*)$ and $h(\Delta)$ such that (\ref{tiaoj}) and
\begin{equation}\label{tj}
h(\Delta)\ge(L^{2q}_{h(\Delta)}\Delta^{\frac{q}{2}})^{-\frac{1}{p-q}}
\end{equation}
holds for any $\Delta\le\Delta_0$, then the continuous-time MTEM methods satisfy
\begin{equation}\label{shou}\mathbb{E}|x(T)-x_\Delta(T)|^q\le C(q,T)L^{2q}_{h(\Delta)}\Delta^{q/2}\ \textrm{and}\ \mathbb{E}|x(T)-\bar{x}_\Delta(T)|^q\le C(q,T)L^{2q}_{h(\Delta)}\Delta^{q/2}.\end{equation}
\end{Theorem}

For the convergence rates over the time interval $[0,T],$ we have to introduce an additional assumption.

Suppose there exist $r\ge2$ and $\bar{K}>0$ such that
\begin{equation}\label{gzeng}
|g(x,y)|^2\le\bar{K}(1+|x|^r+|y|^r), \forall x,y\in \mathbb{R}^d.
\end{equation}

\begin{Theorem}\label{conv1}
Assume that all conditions in Theorem \ref{conv} hold. If (\ref{gzeng}) holds for some $r$ satisfies $2\le r<p-2$ and $2< q\le p-r$ then there exists $C(q,T)$ (independent of $\Delta$) such that
\begin{equation}\label{shou1}\mathbb{E}\sup_{0\le t\le T}|x(t)-x_\Delta(t)|^q\le C(q,T)L^{2q}_{h(\Delta)}\Delta^{q/2}\end{equation}
and if further $2<q<4$, then
\begin{equation}\label{shou2}
\mathbb{E}\sup_{0\le t\le T}|x(t)-\bar{x}_\Delta(t)|^q\le C(q,T)L^{q}_{h(\Delta)}\Delta^{q/2-1}.
\end{equation}
\end{Theorem}

\begin{Remark}
In \cite{GMY}, the authors considered strong convergence of the truncated EM method for SDDEs. However, they only obtained the strong convergence rate at fixed time $T$, while the strong convergence rate over a time interval $[0,T]$ is not considered. Moreover, the weak monotonicity condition is needed (see Assumption 5.1 in \cite{GMY}). In \cite{TY}, the authors obtained the strong convergence rate over a  time interval $[0,T]$, but they also need the weak monotonicity condition A8 (similar to Assumption 5.1 in \cite{GMY}). We only need the local Lipschitz condition (\ref{local}), Khasminskii-type condition (\ref{dandiao}) and (\ref{gzeng}) to make sure the numerical scheme $x_{\Delta}(t)$ strongly converges to $x(t)$ on $[0,T]$ in the sense (\ref{shou1}) with the rate $C(q,T)L^{2q}_{h(\Delta)}\Delta^{q/2}$.
\end{Remark}

\section{Some useful lemmas}

\begin{Lemma}\label{l0}
Suppose the local Lipschitz condition (\ref{local}) holds. Then for any fixed $\Delta>0$,
\begin{equation}\label{global}|f_\Delta(x,y)-f_\Delta(\bar{x},\bar{y})|
\vee|g_\Delta(x,y)-g_\Delta(\bar{x},\bar{y})|\le 5L_{h(\Delta)}(|x-\bar{x}|+|y-\bar{y}|)\end{equation}
for any $x,y,\bar{x},\bar{y}\in\mathbb{R}^d$.
\end{Lemma}

For the proof, see \cite{L} Lemma 3.1.

\begin{Lemma}\label{l1} For $\Delta$ small enough, condition (\ref{dandiao}) implies
\begin{equation}\label{bijin}\aligned 2\langle x-D(y),f_\Delta(x,y)\rangle+(p-1)|g_\Delta(x,y)|^2&\le 2K(1+|x|^2+|y|^2)\endaligned\end{equation}
for any $x,y\in\mathbb{R}^d.$
\end{Lemma}

\textbf{Proof}\ \ On one hand, (\ref{bijin}) holds naturally by (\ref{dandiao}) and the definitions of $f_\Delta$ and $g_\Delta$ if $|x|\vee|y|\le h(\Delta)$.

On the other hand, if $|x|\vee|y|>h(\Delta)$, then
\[\begin{array}{*{20}{l}}
2\langle x-D(y),f_\Delta(x,y)\rangle+|g_\Delta(x,y)|^2 &=2\left\langle x-D(y),\frac{|x|\vee|y|}{h(\Delta)}f\left(\frac{h(\Delta)}{|x|\vee|y|}x,\frac{h(\Delta)}{|x|\vee|y|}y\right)\right\rangle\\
&\quad+\frac{|x|^2\vee|y|^2}{h^2(\Delta)}\left|g\left(\frac{h(\Delta)}{|x|\vee|y|}x,\frac{h(\Delta)}{|x|\vee|y|}y\right)\right|^2\\
&=2\langle x-D(y),\frac{1}{a}f(ax,ay)\rangle+\frac{1}{a^2}|g(ax,ay)|^2\\&
=\frac{1}{a^2}(2\langle ax-aD(\frac{ay}{a}),f(ax,ay)\rangle+|g(ax,ay)|^2)
\end{array}\]
where $a=\frac{h(\Delta)}{|x|\vee|y|}$.
Since $h(\Delta)\ge1$ for sufficiently small $\Delta,$ then by using (\ref{dandiao}), it follows that
\[\begin{array}{*{20}{l}}
2\langle x-D(y),f_\Delta(x,y)\rangle+|g_\Delta(x,y)|^2&\le\frac{1}{a^2}\cdot K(1+|ax|^2+|ay|^2)\\
&=\frac{K_1(|x|^2+|y|^2)}{h^2(\Delta)}+K(|x|^2+|y|^2)\\&
\le2K(1+|x|^2+|y|^2)
\end{array}\]
as required.  $\square$

Now let us state the following two important lemmas. First, we have

\begin{Lemma}\label{ju}
Under conditions (\ref{local}), (\ref{Lip}) and (\ref{dandiao}), the NSDDE (\ref{nsdde}) has a unique global solution $x(t)$ and, moreover, there exists constant $C(p,T)$ (independent of $\Delta$) such that
$$\sup_{0\le t\le T}\mathbb{E}|x(t)|^p\le C(p,T)<\infty,\ \forall T>0.$$ If we define the stopping time
$$\tau_R=\inf\{t\ge0,|x(t)|\ge R\},\ \inf\emptyset=\infty,$$
then for any $T>0,$
$$P(\tau_R\le T)\le\frac{C}{R^p}.$$
\end{Lemma}

\textbf{Proof} By It\^o's formula, we have
$$\aligned\mathbb{E}|x(t)-D(x(t-\tau))|^p&\le \mathbb{E}|\xi(0)-D(\xi(-\tau))|^p+p\mathbb{E}\int_0^t|x(s)-D(x(s-\tau))|^{p-2}\\&
\qquad\qquad\qquad\qquad\qquad\qquad\qquad\times F(x(s),x(s-\tau))ds,
\endaligned$$
where $F(x,y)=\langle x-D(y),f(x,y)\rangle+\frac{p-1}{2}|g(x,y)|^2.$

Then by (\ref{dandiao}), it follows that
$$\aligned&\quad\mathbb{E}|x(t)-D(x(t-\tau))|^p\\&
\le \mathbb{E}|\xi(0)-D(\xi(-\tau))|^p+p K\mathbb{E}\int_0^t|x(s)-D(x(s-\tau))|^{p-2}\\&
\quad\qquad\qquad\qquad\qquad\qquad\qquad\times(1+|x(s)|^2+|x(s-\tau)|^2)ds\\&
\le \mathbb{E}|\xi(0)-D(\xi(-\tau))|^p+p K\mathbb{E}\int_0^t(1+|x(s)-D(x(s-\tau))|^{p})ds\\&
\quad+pK\mathbb{E}\int_0^t\left[\frac{p-2}{p}|x(s)-D(x(s-\tau))|^{p}+|x(s)|^p+|x(s-\tau)|^p\right]ds.
\endaligned$$
We have used Young's inequality in the last inequality.

Thus
$$\aligned\mathbb{E}|x(t)-D(x(t-\tau))|^p&
\le \mathbb{E}|\xi(0)-D(\xi(-\tau))|^p+pKT+pK||\xi||^p_p\\&\quad
+ (2p-2)K\mathbb{E}\int_0^t|x(s)-D(x(s-\tau))|^{p}ds\\&
\quad+2pK\mathbb{E}\int_0^t|x(s)|^pds.\\&
\le \mathbb{E}|\xi(0)-D(\xi(-\tau))|^p+pKT+pK(1+2pu^p)||\xi||^p_p\\&\quad
+[(2p-2)K(1+pu^2)+2pK]\mathbb{E}\int_0^t|x(s)|^pds.
\endaligned$$

Then we have
\begin{equation}\label{sqj}
\aligned \sup_{0\le s\le t}\mathbb{E}|x(s)-D(x(s-\tau))|^p&\le \mathbb{E}|\xi(0)-D(\xi(-\tau))|^p+pKT+pK(1+2pu^p)||\xi||^p_p\\&\quad
+[(2p-2)K(1+pu^2)+2pK]\int_0^t\sup_{0\le r\le s}\mathbb{E}|x(r)|^pds.\endaligned
\end{equation}

On the other hand, for any $c>0$,
\begin{equation}\label{sqj1}\aligned\sup_{0\le s\le t}\mathbb{E}|x(s)|^p&\le (1+c)^{p-1}\sup_{0\le s\le t}\mathbb{E}|x(s)-D(x(s-\tau))|^p\\&
\quad+(\frac{1+c}{c})^{p-1}u^p(||\xi||^p_p+\sup_{0\le s\le t}\mathbb{E}|x(s)|^p).\endaligned\end{equation}

Since $0<u<1$, then we can take $c$ large enough such that $(\frac{1+c}{c})^{p-1}u^p<1$. So
\begin{equation}\label{sqj2}\aligned\sup_{0\le s\le t}\mathbb{E}|x(s)|^p&\le \frac{(1+c)^{p-1}\sup_{0\le s\le t}\mathbb{E}|x(s)-D(x(s-\tau))|^p+(\frac{1+c}{c})^{p-1}u^p||\xi||^p_p}{1-(\frac{1+c}{c})^{p-1}u^p}.\endaligned\end{equation}

Gronwall lemma and (\ref{sqj}) and (\ref{sqj2}) imply that
$$\sup_{0\le t\le T}\mathbb{E}|x(t)|^p\le C(p,T),$$
as required.

Now let us prove
$$P(\tau_R\le T)\le\frac{C}{R^p}.$$

Since
$$\mathbb{E}|x(t\wedge\tau_R)|^p\ge R^pP(\tau_R\le T),$$
so we only need to prove
\begin{equation}\label{sqj5}\sup_{0\le t\le T}\mathbb{E}|x(s\wedge\tau_R)|^p\le C.\end{equation}

As in the above proof, we have
\begin{equation}\label{sqj30}
\aligned \sup_{0\le s\le t}\mathbb{E}|y(s\wedge\tau_R)|^p&\le \mathbb{E}|\xi(0)-D(\xi(-\tau))|^p+pKT+pK(1+2pu^p)||\xi||^p_p\\&\quad
+[(2p-2)K(1+pu^2)+2pK]\int_0^t\sup_{0\le r\le s}\mathbb{E}|x(r\wedge\tau_R)|^pds,\endaligned
\end{equation}
where $y(s\wedge\tau_R)=x(s\wedge\tau_R)-D(x(s\wedge\tau_R-\tau))$
and
\begin{equation}\label{sqj40}\aligned\sup_{0\le s\le t}\mathbb{E}|x(s\wedge\tau_R)|^p&\le \frac{(1+c)^{p-1}\sup_{0\le s\le t}\mathbb{E}|y(s\wedge\tau_R)|^p+(\frac{1+c}{c})^{p-1}u^p||\xi||^p_p}{1-(\frac{1+c}{c})^{p-1}u^p}.\endaligned\end{equation}

Gronwall lemma, (\ref{sqj30}) and (\ref{sqj40}) yield the required (\ref{sqj5}). $\square$

As a similar result of Lemma \ref{ju}, we have the following moment property for the MTEM method (\ref{num2}).

\begin{Lemma}\label{temju}
Assume that (\ref{local}), (\ref{Lip}) and (\ref{dandiao}) hold for $p>2$. Then there exist $0<\Delta_0\le\Delta^*$ and a constant $C(p,T)>0$ (independent of $\Delta$) such that for any $\Delta\in(0,\Delta_0]$, the MTEM method (\ref{num2}) satisfies
\begin{equation}\label{sj}\sup_{0<\Delta\le\Delta_0}\sup_{0\le t\le T}\mathbb{E}|x_\Delta(t)|^p\le C(p,T)<\infty,\ \forall T>0.\end{equation}

Define the stopping time
$$\rho_{\Delta,R}=\inf\{t\ge0,|x_\Delta(t)|\ge R\}.$$ Then for any $R>|x_0|$ and $\Delta\in(0,\Delta^*)$ ($\Delta^*$ small enough), we have
\begin{equation}\label{sj1}P(\rho_{\Delta,R}\le T)\le\frac{C}{R^p}.\end{equation}
\end{Lemma}

\textbf{Proof}\ Let us first prove (\ref{sj}).

Denote $y_\Delta(t)=x_\Delta(t)-D(\bar{x}_\Delta(t-\tau))$, $\bar{y}_\Delta(t)=\bar{x}_\Delta(t)-D(\bar{x}_\Delta(t-\tau))$. By It\^o formula and Lemma \ref{l1}, for any $0\le t\le T,$
$$\aligned&\quad\mathbb{E}(|y_\Delta(t)|^p)\\&\le|y_\Delta(0)|^p+\frac{p}{2}\mathbb{E}\int_0^{t}|y_\Delta(s)|^{p-2}
(2\langle \bar{y}_\Delta(s),f_\Delta(\bar{x}_\Delta(s),\bar{x}_\Delta(s-\tau))\rangle\\&
\qquad\qquad\qquad\qquad\qquad+(p-1)|g_\Delta(\bar{x}_\Delta(s),\bar{x}_\Delta(s-\tau))|^2)ds\\&
\quad+p\mathbb{E}\int_0^{t}|y_\Delta(s)|^{p-2}
\langle x_\Delta(s)-\bar{x}_\Delta(s),f_\Delta(\bar{x}_\Delta(s),\bar{x}_\Delta(s-\tau))\rangle ds\\&
\le|y_\Delta(0)|^p+\frac{p}{2}\mathbb{E}\int_0^{t}|y_\Delta(s)|^{p-2}\cdot 2K(1+|\bar{x}_\Delta(s)|^2+|\bar{x}_\Delta(s-\tau)|^2)ds\\&
\quad+p\mathbb{E}\int_0^{t}|y_\Delta(s)|^{p-2}
|x_\Delta(s)-\bar{x}_\Delta(s)|\cdot (5L_{h(\Delta)}(|\bar{x}_\Delta(s)|+|\bar{x}_\Delta(s-\tau)|)+|f(0,0)|)ds\\&
\le|x_0|^p+C_p\mathbb{E}\int_0^{t}(1+|y_\Delta(s)|^{p}+|\bar{x}_\Delta(s)|^p)ds\\&
\quad+C_p\mathbb{E}\int_0^{t}(|y_\Delta(s)|^{p}+L_{h(\Delta)}^p|x_\Delta(s)-\bar{x}_\Delta(s)|^p+
|\bar{x}_\Delta(s)|^p+|\bar{x}_\Delta(s-\tau)|^p)ds
\\&\quad+C_p|f(0,0)|\mathbb{E}\int_0^{t}(1+|y_\Delta(s)|^{p}+|x_\Delta(s)-\bar{x}_\Delta(s)|^p)ds\\&
\le C_p+C_p\mathbb{E}\int_0^{t}|y_\Delta(s)|^{p}ds+C_p\mathbb{E}\int_0^{t}\sup_{0\le r\le s}\mathbb{E}|\bar{x}_\Delta(r)|^pds\\&
\quad+(C_p L_{h(\Delta)}^p+1)\int_0^{t}\sup_{0\le r\le s}\mathbb{E}|x_\Delta(r)-\bar{x}_\Delta(r)|^pds.\endaligned$$

Notice that for any $0\le s\le t,$ there exists $k\le[\frac{t}{\Delta}]$ such that $k\Delta\le s<(k+1)\Delta.$ Thus
$$x_\Delta(s)-\bar{x}_\Delta(s)=x_\Delta(s)-X_k=f_\Delta(X_k,X_{k-m})(s-k\Delta)+g_\Delta(X_k,X_{k-m})(B(s)-B(k\Delta)).$$

So we have
$$\aligned\mathbb{E}|x_\Delta(s)-\bar{x}_\Delta(s)|^p&\le C_p\left[\Delta^p\mathbb{E}(|f_\Delta(X_k,X_{k-m})|^p)\right.\\&
\qquad\left.+\mathbb{E}(|g_\Delta(X_k,X_{k-m})|^p)\mathbb{E}(|B(s)-B(k\Delta)|^p|\mathscr{F}_{k\Delta})\right].\endaligned$$

Since $f_\Delta$ and $g_\Delta$ satisfy the global Lipschitz condition (\ref{global}), and notice that $B(t)-B(k\Delta)$ is independent of $\mathscr{F}_{k\Delta},$ then
$$\aligned\mathbb{E}|x_\Delta(s)-\bar{x}_\Delta(s)|^p&\le C_p\big(\Delta^p\mathbb{E}(5L_{h(\Delta)}(|X_k|+|X_{k-m}|)+|f(0,0)|)^p\\&\quad
+\mathbb{E}(5L_{h(\Delta)}(|X_k|+|X_{k-m}|)+|g(0,0)|)^p\Delta^\frac{p}{2}\big)\\&
\le C_pL_{h(\Delta)}^p\Delta^p(\mathbb{E}(|X_k|^p)+\mathbb{E}(|X_{k-m}|^p))+C_p\Delta^p|f(0,0)|^p\\&\quad
+C_pL_{h(\Delta)}^p\Delta^\frac{p}{2}(\mathbb{E}(|X_k|^p)+\mathbb{E}(|X_{k-m}|^p))+C_p\Delta^\frac{p}{2}|g(0,0)|^p.\endaligned$$

Therefore, for any $t\le T,$
\begin{equation}\label{k2}\aligned\sup_{0\le s\le t}\mathbb{E}|x_\Delta(s)-\bar{x}_\Delta(s)|^p&\le C_p(L^p_{h(\Delta)}\Delta^\frac{p}{2}\sup_{0\le k\le[\frac{t}{\Delta}]}\mathbb{E}(|X_k|^p)+\Delta^\frac{p}{2}),\endaligned\end{equation}
where $C_p$ is a positive constant (independent of $\Delta$) which might change values from line to line.

Since $L^4_{h(\Delta)}\Delta\to0$ as $\Delta\to0$, then $(C_p L_{h(\Delta)}^p+1)L^p_{h(\Delta)}\Delta^\frac{p}{2}$ is bounded for $\Delta\in(0,\Delta_0]$, therefore we have
$$\aligned\quad\mathbb{E}(|y_\Delta(t)|^p)&\le
 C_p+C_p\mathbb{E}\int_0^{t}|y_\Delta(s)|^{p}ds+C_p\mathbb{E}\int_0^{t}\sup_{0\le r\le s}\mathbb{E}|\bar{x}_\Delta(r)|^pds\\&
\quad+(C_p L_{h(\Delta)}^p+1)[C_pL^p_{h(\Delta)}\Delta^\frac{p}{2}\int_0^{t}\sup_{0\le s\le t}\mathbb{E}(|x_\Delta(s)|^p)ds+C_pT\Delta^\frac{p}{2}]\\&
\le C_p+C_p\int_0^{t}\sup_{0\le s\le t}\mathbb{E}(|x_\Delta(s)|^p)ds.\endaligned$$

On the other hand, similar to (\ref{sqj2}), we can take $c$ large enough such that $(\frac{1+c}{c})^{p-1}u^p<1$. So
\begin{equation}\label{sqj33}\aligned\sup_{0\le s\le t}\mathbb{E}|x_\Delta(s)|^p&\le \frac{(1+c)^{p-1}\sup_{0\le s\le t}\mathbb{E}|y_\Delta(s)|^p+(\frac{1+c}{c})^{p-1}u^p||\xi||^p_p}{1-(\frac{1+c}{c})^{p-1}u^p}.\endaligned\end{equation}

Thus
$$\aligned\quad\sup_{0\le s\le t}\mathbb{E}|x_\Delta(s)|^p&\le
 C_p+C_p\mathbb{E}\int_0^{t}\sup_{0\le r\le s}\mathbb{E}|x_\Delta(r)|^pds.\endaligned$$

Gronwall inequality yields (\ref{sj}).

Now let us prove (\ref{sj1}).

Let $y_\Delta(t)$, $\bar{y}_\Delta(t)$ are defined as above and $\rho_{\Delta,R}=\rho.$ By It\^o formula and Lemma \ref{l1}, for any $0\le t\le T,$

$$\aligned\mathbb{E}(|y_\Delta(t\wedge\rho)|^p)&
\le|x_0|^p+C_p\mathbb{E}\int_0^{t\wedge\rho}(1+|y_\Delta(s)|^{p}+|\bar{x}_\Delta(s)|^p)ds\\&
\quad+C_p\mathbb{E}\int_0^{t\wedge\rho}(|y_\Delta(s)|^{p}+L_{h(\Delta)}^p|x_\Delta(s)-\bar{x}_\Delta(s)|^p+
|\bar{x}_\Delta(s)|^p+|\bar{x}_\Delta(s-\tau)|^p)ds
\\&\quad+C_p|f(0,0)|\mathbb{E}\int_0^{t\wedge\rho}(1+|y_\Delta(s)|^{p}+|x_\Delta(s)-\bar{x}_\Delta(s)|^p)ds\\&
\le C_p+C_p\mathbb{E}\int_0^{t\wedge\rho}|y_\Delta(s)|^{p}ds+C_pT\sup_{0\le s\le T}\mathbb{E}|\bar{x}_\Delta(s)|^p\\&
\quad+(C_p T L_{h(\Delta)}^p+1)\sup_{0\le s\le T}\mathbb{E}|x_\Delta(s)-\bar{x}_\Delta(s)|^p.\endaligned$$

Then by (\ref{k2}) and (\ref{sj}), we have

$$\mathbb{E}(|y_\Delta(t\wedge\rho)|^p)\le C_p+C_p\int_0^{t}\mathbb{E}|y_\Delta(s\wedge\rho)|^pds+(C_p T L_{h(\Delta)}^p+1)L_{h(\Delta)}^p\Delta^\frac{p}{2}.$$

Gronwall's lemma yields that
$$\sup_{0\le s\le T}\mathbb{E}(|y_\Delta(s\wedge\rho)|^p)\le C(p,T)<\infty.$$

Then similar to (\ref{sqj2}), we can take $c$ large enough such that $(\frac{1+c}{c})^{p-1}u^p<1$. So
\begin{equation}\label{sqj4}\aligned\sup_{0\le s\le t}\mathbb{E}|x_\Delta(s\wedge\rho)|^p&\le \frac{(1+c)^{p-1}\sup_{0\le s\le t}\mathbb{E}|y_\Delta(s\wedge\rho)|^p+(\frac{1+c}{c})^{p-1}u^p||\xi||^p_p}{1-(\frac{1+c}{c})^{p-1}u^p}.\endaligned\end{equation}

This implies the required assertion easily. $\square$

\section{Convergence rate at fixed time $T$}

Let us first present a lemma which will play a key role in the proof of Theorem \ref{conv}.

\begin{Lemma}\label{jubu}
Suppose (\ref{local}) and (\ref{Lip}) hold for $p>2,$ and for any $2<q<p,$
$$\mathbb{E}\sup_{-m\le k\le -1}\sup_{k\Delta\le s\le (k+1)\Delta}|\xi(s)-\xi(k\Delta)|^q\le \hat{K}\Delta^\frac{q}{2}.$$
Set $$\theta_{\Delta,R}=\tau_R\wedge\rho_{\Delta,R}\quad and\quad e_\Delta(t)=x(t)-x_\Delta(t)\ for\ t\ge0.$$
Then for any $\Delta\in(0,\Delta^*)$ and any $R\le h(\Delta^*)$, there exists $C(q,T)>0$ (independent of $\Delta$) such that
$$\sup_{0\le t\le T}\mathbb{E}(|e_\Delta(t\wedge\theta_{\Delta,R})|^q)\le C(q,T)L^{2q}_{h(\Delta)}\Delta^{\frac{q}{2}}.$$
\end{Lemma}

\textbf{Proof}\ Define the truncated functions
$$F_R(x,y)=f_{h^{-1}(R)}(x,y)\ \textrm{and}\ G_R(x,y)=g_{h^{-1}(R)}(x,y), \forall x,y\in\mathbb{R}^d,$$
where $f_{h^{-1}(R)}$ is defined in (\ref{dy}) with $\Delta$ replaced by $h^{-1}(R).$ By Lemma (3.1), $F_R$ and $G_R$ are globally Lipschitz continuous for any fixed $R$ ($\ge L^{-1}$(1)), where $L^{-1}$ is the inverse function of $L_R$ when it is seen as a function of $R.$

Without loss of generality, suppose $\Delta^*$ is sufficiently small such that
$$h(\Delta^*)=L^{-1}(L_R\exp\{\frac{2^{q-1}(T^\frac{1}{q-1}+4)}{q}\})\ge R.$$

Then for those $x,y\in\mathbb{R}^d$ with $|x|\vee|y|\le R$ and all $\Delta\in(0,\Delta^*]$, we have
$$F_R(x,y)=f_{h^{-1}(R)}(x,y)=f(x,y)=f_\Delta(x,y).$$

Similarly, we have $$G_R(x,y)=g_\Delta(x,y).$$

Now consider NSDDE
\begin{equation}\label{fz}
d[z(t)-D(z(t-\tau))
=F_R(z(t),z(t-\tau))dt+G_R(z(t),z(t-\tau))dB_t,t\ge0
\end{equation}
with $z(\theta)=\xi(\theta)$ on $\theta\in[-\tau,0].$ Since $F_R$ and $G_R$ are globally Lipschitz continuous (with Lipschitz constant $5L_R$) for any fixed $R$, then (\ref{fz}) has a unique global solution $z(t)$ on $t\ge\tau.$ Thus
\begin{equation}\label{weiyi}P(x(t\wedge\tau_R)=z(t\wedge\tau_R), \forall t\in[0,T])=1.\end{equation}

On the other hand, similar to (\ref{num1}) and (\ref{num2}), we can define $\bar{z}_\Delta(t), z_\Delta(t)$ in the same way for NSDDE (\ref{fz}). We also have
\begin{equation}\label{weiyi1}P(x_\Delta(t\wedge\tau_R)=z_\Delta(t\wedge\tau_R), \forall t\in[0,T])=1.\end{equation}

We claim that
\begin{equation}\label{zhengti}\mathbb{E}\sup_{0\le t\le T}|z(t)-z_\Delta(t)|^q\le C(q,T)L^{2q}_{h(\Delta)}\Delta^\frac{q}{2}.\end{equation}

Let $y(t)=z(t)-D(z(t-\tau)), y_\Delta(t)=z_\Delta(t)-D(\bar{z}_\Delta(t-\tau))$. Then for any $c,c'>0$
$$\aligned|z(t)-z_\Delta(t)|^q&\le (1+c)^{q-1}|y(t)-y_\Delta(t)|^q+(\frac{1+c}{c})^{q-1}u^q|z(t-\tau)-\bar{z}_\Delta(t-\tau)|^q\\&
\le (1+c)^{q-1}|y(t)-y_\Delta(t)|^q+(\frac{(1+c)(1+c')}{c})^{q-1}u^q|z(t-\tau)-z_\Delta(t-\tau)|^q\\&
\quad+(\frac{(1+c)(1+c')}{cc'})^{q-1}u^q|z_\Delta(t-\tau)-\bar{z}_\Delta(t-\tau)|^q.\endaligned$$

Choose $c$ sufficiently large and $c'$ sufficiently small such that $c_0:=(\frac{(1+c)(1+c')}{c})^{q-1}u^q<1$, and denote $c_1=(\frac{(1+c)(1+c')}{cc'})^{q-1}u^q$. Then we have
$$\aligned\mathbb{E}\sup_{0\le s\le t}|z(s)-z_\Delta(s)|^q&\le (1+c)^{q-1}\mathbb{E}\sup_{0\le s\le t}|y(s)-y_\Delta(s)|^q+c_0\mathbb{E}\sup_{0\le s\le t}|z(s-\tau)-z_\Delta(s-\tau)|^q\\&
\quad+c_1\mathbb{E}\sup_{0\le s\le t}|z(s-\tau)-\bar{z}_\Delta(s-\tau)|^q\\&
\le (1+c)^{q-1}\mathbb{E}\sup_{0\le s\le t}|y(s)-y_\Delta(s)|^q+c_0\mathbb{E}\sup_{0\le s\le t}|z(s)-z_\Delta(s)|^q\\&
\quad+c_1\mathbb{E}\sup_{0\le s\le t}|z_\Delta(s)-\bar{z}_\Delta(s)|^q\\&
\quad+c_1\mathbb{E}\sup_{-m\le k\le -1}\sup_{k\Delta\le s\le (k+1)\Delta}|\xi(s)-\xi(k\Delta)|^q.\endaligned$$

So
\begin{equation}\label{kz0}\aligned\mathbb{E}\sup_{0\le s\le t}|z(s)-z_\Delta(s)|^q&\le \frac{(1+c)^{q-1}}{1-c_0}\mathbb{E}\sup_{0\le s\le t}|y(s)-y_\Delta(s)|^q\\&
\quad+\frac{c_1}{1-c_0}\mathbb{E}\sup_{0\le s\le t}|z_\Delta(s)-\bar{z}_\Delta(s)|^q+\frac{c_1}{1-c_0}\hat{K}\Delta^{\frac{q}{2}}.\endaligned\end{equation}

As in (\ref{k2}), we have
\begin{equation}\label{k4}\sup_{0\le t\le T}\mathbb{E}|z_\Delta(t)-\bar{z}_\Delta(t)|^p\le C(p,T)L^p_{h(\Delta)}\Delta^\frac{p}{2},\end{equation}

Now by It\^o's formula, H\"{o}lder's inequality and BDG inequality, it follows that for $0\le t\le T$
$$\aligned&\quad\mathbb{E}\sup_{0\le s\le t}|y(s)-y_\Delta(s)|^q\\&\le 2^{q-1}T^{\frac{1}{q-1}}\mathbb{E}\int_0^t|F_R(z(s),z(s-\tau))-F_R(\bar{z}_\Delta(s),\bar{z}_\Delta(s-\tau))|^qds\\&
\quad+2^{q-1}\mathbb{E}\int_0^t|G_R(z(s),z(s-\tau))-G_R(\bar{z}_\Delta(s),\bar{z}_\Delta(s-\tau))|^qds\\&
\le 2^{q-1}(T^{\frac{1}{q-1}}+4)L_R^q\int_0^t(\mathbb{E}|z(s)-\bar{z}_\Delta(s)|^q+\mathbb{E}|z(s-\tau)-\bar{z}_\Delta(s-\tau)|^q)ds\\&
\le 2^{q}(T^{\frac{1}{q-1}}+4)L_R^q\int_0^t\mathbb{E}|z(s)-z_\Delta(s)|^qds+C(q,T)L_R^q\sup_{0\le t\le T}\mathbb{E}|z_\Delta(t)-\bar{z}_\Delta(t)|^q\\&
\quad+C(q,T)L_R^q\sum_{k=-m}^{-1}\int_{k\Delta}^{(k+1)\Delta}\mathbb{E}|\xi(s)-\xi(k\Delta)|^qds\\&
\le 2^{q}(T^{\frac{1}{q-1}}+4)L_R^q\int_0^t\mathbb{E}|z(s)-z_\Delta(s)|^qds+C(q,T)L_R^{2q}\Delta^{\frac{q}{2}}\\&
\quad+C(q,T)L_R^q\hat{K}\Delta^{\frac{q}{2}}.
\endaligned$$

Then (\ref{kz0}) and Gronwall lemma yields
\begin{equation}\label{zhengti1}\mathbb{E}\sup_{0\le t\le T}|z(t)-z_\Delta(t)|^q\le H(R,T,\xi)\Delta^\frac{q}{2},\end{equation}
where
$$\aligned H(R,T,\xi)&:=(C(q,T)L_R^{2q}+C(q,T)L_R^q\hat{K})e^{2^{q}(T^{\frac{1}{q-1}}+4)TL_R^q}\\&
\le C(q,T) L_R^{2q}\exp\{2^{q}(T^{\frac{1}{q-1}}+4)TL_R^q\}\\&
=C(q,T)L^{2q}_{h(\Delta^*)}\le C(q,T)L^{2q}_{h(\Delta)}.\endaligned$$

Hence (\ref{weiyi}), (\ref{weiyi1}) and (\ref{zhengti}) implies
\begin{equation}\label{kz01}\aligned\mathbb{E}\sup_{0\le s\le t}|x(s\wedge\theta_{\Delta,R})-x_\Delta(s\wedge\theta_{\Delta,R})|^q&\le C(q,T)L^{2q}_{h(\Delta)}\Delta^{\frac{q}{2}},\endaligned\end{equation}

This completes the proof. $\square$

Now we are ready to prove Theorem \ref{conv}.

\textbf{Proof of Theorem \ref{conv}}\ Let $\tau_R, \rho_{\Delta,R}, \theta_{\Delta,R}$ and $e_\Delta(t)$ be the same as before.
Then by Young's inequality, we have that for any $\delta>0,$
\begin{equation}\label{bds}\aligned\mathbb{E}(|e_\Delta(T)|^q)&\le \mathbb{E}(|e_\Delta(T)|^q1_{\{\theta_{\Delta,R}>T\}})+ \mathbb{E}(|e_\Delta(T)|^q1_{\{\theta_{\Delta,R}\le T\}})\\&\le\mathbb{E}(|e_\Delta(T)|^q1_{\{\theta_{\Delta,R}>T\}})+ \frac{q\delta}{p}\mathbb{E}(|e_\Delta(T)|^p)+\frac{p-q}{p\delta^{q/(p-q)}}P(\theta_{\Delta,R}\le T)\\&
\le\mathbb{E}(|e_\Delta(T\wedge\theta)|^q)+ \frac{q\delta C}{p}\left(\mathbb{E}(|x_\Delta(T)|^p)+\mathbb{E}(|x(T)|^p)\right)+\frac{p-q}{p\delta^{q/(p-q)}}P(\theta_{\Delta,R}\le T)\endaligned$$
where $C$ is a positive constant (independent of $\Delta$) which might change the value from line to line. We have used the fact that
$$\mathbb{E}(|e_\Delta(T)|^p)\le C(\mathbb{E}(|x_\Delta(T)|^p)+\mathbb{E}(|x(T)|^p))\end{equation}
in the last inequality.

By Lemma \ref{ju} and \ref{temju}, we have
$$\mathbb{E}(|x_\Delta(T)|^p)+\mathbb{E}(|x(T)|^p)\le C,$$
and
$$P(\theta_{\Delta,R}\le T)\le P(\tau_R\le T)+P(\rho_{\Delta,R}\le T)\le\frac{C}{R^p}.$$

Thus,
$$\mathbb{E}(|e_\Delta(T)|^q)\le\mathbb{E}(|e_\Delta(T\wedge\theta)|^q)+\frac{qC\delta}{p}+\frac{C(p-q)}{pR^p\delta^{q/(p-q)}}$$
holds for any $\Delta\in(0,\Delta^*), R>|x_0|$ and $\delta>0$. Then we can choose $\delta=L^{2q}_{h(\Delta)}\Delta^\frac{q}{2}$ and $R=(L^{2q}_{h(\Delta)}\Delta^\frac{q}{2})^{-\frac{1}{p-q}}$ to get
$$\mathbb{E}(|e_\Delta(T)|^q)\le\mathbb{E}(|e_\Delta(T\wedge\theta)|^q)+CL^{2q}_{h(\Delta)}\Delta^\frac{q}{2}.$$

But by condition (\ref{tj}), we have
$$h(\Delta)\ge (L^{2q}_{h(\Delta)}\Delta^\frac{q}{2})^{-\frac{1}{p-q}}=R.$$

Then by Lemma \ref{jubu},
$$\mathbb{E}(|e_\Delta(T)|^q)\le C(q,T)L^{2q}_{h(\Delta)}\Delta^\frac{q}{2},$$
where $C$ is a positive constant depends on $q$ and $T$. This is the first inequality of (\ref{shou}).

For the second inequality, since $q<p,$ by H\"older inequality, it follows easily from the above inequality and (\ref{k2}) in Lemma \ref{temju} that
$$\aligned\mathbb{E}(|x(T)-\bar{x}_\Delta(T)|^q)&\le C_q\left(\mathbb{E}(|e_\Delta(T)|^q)+\mathbb{E}(|x_\Delta(T)-\bar{x}_\Delta(T)|^q)\right)\\&\le C_q\left(\mathbb{E}(|e_\Delta(T)|^q)+[\sup_{0\le t\le T}\mathbb{E}(|x_\Delta(t)-\bar{x}_\Delta(t)|^p)]^\frac{q}{p}\right)\\&
\le C(q,T)L^{2q}_{h(\Delta)}\Delta^\frac{q}{2}+(C(q,T)L^{p}_{h(\Delta)}\Delta^\frac{p}{2})^\frac{q}{p}\le C(q,T)L^{2q}_{h(\Delta)}\Delta^\frac{q}{2}.\endaligned$$
We complete the proof. $\square$

\section{Convergence rates over the time interval $[0,T]$}

First of all, let us prove a similar Lemma to Lemma \ref{ju}.

\begin{Lemma}\label{zuida}
Let (\ref{local}), (\ref{dandiao}) and (\ref{gzeng}) hold for $p>2$ and $2<r<p.$ Set $\bar{p}=p-r+2.$ Then
$$\mathbb{E}(\sup_{0\le t\le T}|x(t)|^{\bar{p}})\le C,\ \forall T>0.$$
\end{Lemma}

\textbf{Proof}\ Since for sufficiently large $c>0$,
$$\mathbb{E}(\sup_{0\le t\le T}|x(t)|^{\bar{p}})\le \frac{(1+c)^{\bar{p}-1}\mathbb{E}\sup_{0\le t\le T}|y(t)|^{\bar{p}}+(\frac{1+c}{c})^{\bar{p}-1}u^{\bar{p}}||\xi||^{\bar{p}}}{1-(\frac{1+c}{c})^{\bar{p}-1}u^{\bar{p}}},$$
then we only need to prove
$$\mathbb{E}\sup_{0\le t\le T}|y(t)|^{\bar{p}}\le C(\bar{p},T).$$

Indeed, It\^o's formula and (\ref{dandiao}) imply that
$$\aligned\mathbb{E}(\sup_{0\le t\le T}|y(t)|^{\bar{p}})&\le |y(0)|^{\bar{p}}+\bar{p}K\mathbb{E}\int_0^T|y(s)|^{\bar{p}-2}(1+|x(s)|^2+|x(s-\tau)|^2)ds\\&
\quad+\mathbb{E}(\sup_{0\le t\le T}|M(t)|),\endaligned$$
where $y(t)=x(t)-D(x(t-\tau)),$ and
$$M(t)=\bar{p}\int_0^t|y(s)|^{\bar{p}-2}y^T(s)g(x(s),x(s-\tau))dB(s)$$
is a local martingale with $M(0)=0.$

It is obvious that
$$\mathbb{E}\int_0^T|y(s)|^{\bar{p}-2}(1+|x(s)|^2+|x(s-\tau)|^2)ds\le C(\bar{p},T)<\infty.$$

On the other hand, by Burkholder-Davis-Gundy (BDG) inequality (see e.g. \cite{Ikeda}), it follows that
$$\aligned\mathbb{E}(\sup_{0\le t\le T}|M(t)|)&\le C'\mathbb{E}\left(\int_0^T|y(s)|^{2\bar{p}-2}|g(x(s),x(s-\tau))|^2ds\right)^\frac{1}{2}\\&
\le C'\mathbb{E}\left(\sup_{0\le t\le T}|y(t)|^{\bar{p}}\int_0^T|y(s)|^{\bar{p}-2}\cdot \bar{K}(1+|x(s)|^r+|x(s-\tau)|^r)ds\right)^\frac{1}{2}\\&
\le \frac{1}{2}\mathbb{E}\sup_{0\le t\le T}|y(t)|^{\bar{p}}+\frac{C'^2}{2}\mathbb{E}\int_0^T|y(s)|^{\bar{p}-2}\cdot \bar{K}(1+|x(s)|^r+|x(s-\tau)|^r)ds.\endaligned$$

As in the proof of Lemma \ref{ju}, we have
$$\frac{C'^2}{2}\mathbb{E}\int_0^T|y(s)|^{\bar{p}-2}\cdot \bar{K}(1+|x(s)|^r+|x(s-\tau)|^r)ds\le C''.$$

Thus
$$\aligned\mathbb{E}(\sup_{0\le t\le T}|y(t)|^{\bar{p}})&\le C(\bar{p},T),\endaligned$$

This completes the proof. $\square$

For the discontinuous and continuous-time MTEM methods (\ref{num}) and (\ref{num2}), we have
\begin{Lemma}\label{zuida1}
Let (\ref{local}), (\ref{tiaoj}), (\ref{dandiao}) and (\ref{gzeng}) hold for $p>2$ and $p>r\ge2.$ Set $\bar{p}=p+2-r.$ Then
\begin{equation}\label{zuida2}\sup_{0<\Delta\le\Delta^*}\mathbb{E}(\sup_{0\le t\le T}|x_\Delta(t)|^{\bar{p}})\le C,\ \forall T>0,\end{equation}
and therefore,
\begin{equation}\label{zuida3}\sup_{0<\Delta\le\Delta^*}\mathbb{E}(\sup_{0\le k\le [\frac{T}{\Delta}]}|X_k^\Delta|^{\bar{p}})\le C,\ \forall T>0,\end{equation}
\end{Lemma}

\textbf{Proof}\ Similar to the proof of Lemma \ref{zuida}, we only need to prove
$$\sup_{0<\Delta\le\Delta^*}\mathbb{E}(\sup_{0\le t\le T}|y_\Delta(t)|^{\bar{p}})\le C(\bar{p},T),$$
where $y_\Delta(t)=x_\Delta(t)-D(\bar{x}_\Delta(t-\tau)).$

For any $\Delta\in(0,\Delta^*],$ by It\^o formula and Lemma \ref{l1}, we have
$$\aligned\sup_{0\le t\le T}|y_\Delta(t)|^{\bar{p}}&\le|y_0|^{\bar{p}}+2K\bar{p}\int_0^T|y_\Delta(s)|^{\bar{p}-2}
(1+|\bar{x}_\Delta(s)|^2+\bar{x}_\Delta(s-\tau)|^2)ds\\&
\quad+\bar{p}\int_0^T|y_\Delta(s)|^{\bar{p}-2}|x_\Delta(s)-\bar{x}_\Delta(s)||f_\Delta(\bar{x}_\Delta(s),\bar{x}_\Delta(s-\tau))|ds\\&
\quad+\bar{p}\sup_{0\le t\le T}\left|\int_0^t|y_\Delta(s)|^{\bar{p}-2}\langle y_\Delta(s),g_\Delta(\bar{x}_\Delta(s),\bar{x}_\Delta(s-\tau))dB(s)\rangle\right|.\endaligned$$

Young's inequality and Lemma \ref{temju} imply that
$$\int_0^T|y_\Delta(s)|^{\bar{p}-2}
(1+|\bar{x}_\Delta(s)|^2+\bar{x}_\Delta(s-\tau)|^2)ds<C(\bar{p},T)$$

Moreover, since for $\Delta>0$ small enough,
$$|f_\Delta(x,y)|\le 5L_{h(\Delta)}(|x|+|y|)+|f(0,0)|,$$
then as in the proof of Lemma \ref{temju}, we have
$$\mathbb{E}\int_0^T|y_\Delta(s)|^{\bar{p}-2}|x_\Delta(s)-\bar{x}_\Delta(s)||f_\Delta(\bar{x}_\Delta(s),\bar{x}_\Delta(s-\tau))|ds\le C.$$

So by BDG inequality again and (\ref{gzeng}), we have
$$\aligned&\quad\mathbb{E}\sup_{0\le t\le T}|y_\Delta(t)|^{\bar{p}}\\&
\le C+C(\bar{p},\bar{K})\mathbb{E}\left(\int_0^T|y_\Delta(s)|^{2\bar{p}-2}(1+|\bar{x}_\Delta(s)|^r+|\bar{x}_\Delta(s-\tau)|^r)ds\right)^\frac{1}{2}\\&
\le C+\mathbb{E}\left|\sup_{0\le t\le T}|y_\Delta(t)|^{\bar{p}}\cdot C(\bar{p},\bar{K})\left(\int_0^T|y_\Delta(s)|^{\bar{p}-2}
(1+|\bar{x}_\Delta(s)|^r+|\bar{x}_\Delta(s-\tau)|^r)ds\right)\right|^\frac{1}{2}\\&
\le C+\frac{1}{2}\mathbb{E}\sup_{0\le t\le T}|y_\Delta(t)|^{\bar{p}}+\frac{C^2(\bar{p},\bar{K})}{2}\int_0^T|y_\Delta(s)|^{\bar{p}-2}
(1+|\bar{x}_\Delta(s)|^r+|\bar{x}_\Delta(s-\tau)|^r)ds,
\endaligned$$
where $C$ is a constant (independent of $\Delta$).

Then we have
$$\aligned\mathbb{E}\sup_{0\le t\le T}|y_\Delta(t)|^{\bar{p}}
&\le C+C^2(\bar{p},\bar{K})\mathbb{E}\int_0^T|y_\Delta(s)|^{\bar{p}-2}
(1+|\bar{x}_\Delta(s)|^r+|\bar{x}_\Delta(s-\tau)|^r)ds\le C.\endaligned$$

Then the required assertion (\ref{zuida2}) follows. $\square$

\begin{Lemma}\label{close1} Assume that (\ref{local}), (\ref{tiaoj}), (\ref{dandiao}) and (\ref{gzeng}) hold for $p>2$ and $2\le r< p.$ If $q\le p+2-r,$ then for any $\Delta\in(0,\Delta^*)$, there exists $C>0$ (independent of $\Delta$) such that
\begin{equation}\mathbb{E}\left(\sup_{0\le t\le T}|x_\Delta(t)-\bar{x}_\Delta(t)|^q\right)\le CL_{h(\Delta)}^q\Delta^{\frac{q}{2}-1}.\end{equation}
\end{Lemma}
The proof is almost the same as that of Lemma 5.5 in \cite{LX}, so we omit it here.

Now let us prove Theorem \ref{conv1}.

\textbf{Proof of Theorem \ref{conv1}}\ Let $\theta_{\Delta,R}$ and $e_\Delta(t)$ be the same as before.
As in the proof of Theorem \ref{conv},  by Young's inequality, we have that for any $\delta>0,$
$$\aligned\mathbb{E}(\sup_{0\le t\le T}|e_\Delta(t)|^q)&\le \mathbb{E}(1_{\{\theta_{\Delta,R}>T\}}\sup_{0\le t\le T}|e_\Delta(t)|^q)+ \frac{q\delta}{p}\mathbb{E}(\sup_{0\le t\le T}|e_\Delta(t)|^p)\\&\quad+\frac{p-q}{p\delta^{q/(p-q)}}P(\theta_{\Delta,R}\le T).\endaligned$$

By Lemma \ref{zuida}, \ref{zuida1},
$$\mathbb{E}(\sup_{0\le t\le T}|e_\Delta(t)|^p)\le C(\mathbb{E}(\sup_{0\le t\le T}|x(t)|^p)+\mathbb{E}(\sup_{0\le t\le T}|x_\Delta(t)|^p))\le C.$$

Then similar to the proof of Theorem \ref{conv}, we have
$$\aligned\mathbb{E}(\sup_{0\le t\le T}|e_\Delta(t)|^q)&\le \mathbb{E}(\sup_{0\le t\le T}|e_\Delta(t\wedge\theta_{\Delta,R})|^q)+ \frac{Cq\delta}{p}+\frac{C(p-q)}{pR^p\delta^{q/(p-q)}}\endaligned$$
holds for any $\Delta\in(0,\Delta^*), \delta>0$ and $R>|x_0|.$

Since we have proved Lemma \ref{jubu}, repeat the proof of Theorem \ref{conv}, we have
$$\mathbb{E}(\sup_{0\le t\le T}|e_\Delta(t)|^q)\le CL^{2q}_{h(\Delta)}\Delta^\frac{q}{2}$$
for $\delta=L^{2q}_{h(\Delta)}\Delta^\frac{q}{2},\quad R=(L^{2q}_{h(\Delta)}\Delta^\frac{q}{2})^{-\frac{1}{p-q}}\le h(\Delta),$ as required. $\square$

\section{Examples}

Now let us present two examples to illustrate our theory.

\textbf{Example 1}
Let $d=1, \tau=1.$ Consider the following scalar NSDDE:
 \begin{equation}\label{sde1}\aligned d[x(t)-\frac{1}{2}\sin x(t-1)]&=(2x(t)-x(t)e^{3x(t)}+\frac{1}{2}\sin
x(t-1))dt\\&\quad+\sqrt{\frac{1}{5}x^2(t)e^{3x(t)}+x^2(t-1)+1}dB_t.\endaligned\end{equation}

Here $f(x,y)=2x-xe^{3x}+\frac{1}{2}\sin y$, $g(x,y)=\sqrt{\frac{1}{5}x^2e^{3x}+y^2+1}$ and $D(y)=\frac{1}{2}\sin y.$ Then neither $f$ nor $g$ is polynomial growing (although both are local Lipschitz continuous).

Moreover, for any $a\in(0,1]$, we have
\begin{equation}
\aligned2\langle
x-aD(\frac{1}{a}y),f(x,y)\rangle+5|g(x,y)|^2&=4x^2-2x^2e^{3x}+x\sin y-2ax\sin \frac{y}{a}\\&\quad+axe^{3x}\sin \frac{y}{a}-\frac{a}{2}\sin y\sin\frac{y}{a}+x^2e^{3x}+5y^2+5\\&
\le4x^2-x^2e^{3x}+|x|+2|x|\\&\quad+|x|e^{3x}+\frac{1}{2}+5y^2+5\\&
\le4x^2+3\cdot\frac{1+x^2}{2}+\frac{11}{2}+5y^2+e^{3x}(|x|-x^2).\endaligned
\end{equation}

Notice that if $x\le0$, we have
$$e^{3x}(|x|-x^2)=e^{3x}[\frac{1}{4}-(|x|-\frac{1}{2})^2]\le \frac{1}{4}e^{3x}\le \frac{1}{4}.$$

If $x>0$, then
$$e^{3x}(|x|-x^2)\le\sup_{0\le x\le1}e^{3x}(x-x^2)=e^{3x}(x-x^2)|_{x=\frac{\sqrt{5}-1}{2}}=(\sqrt{5}-2)e^\frac{3\sqrt{5}-3}{2}<e^2.$$

Thus
\begin{equation}
\aligned2\langle
x-aD(\frac{1}{a}y),f(x,y)\rangle+3|g(x,y)|^2&
\le4x^2+3\cdot\frac{1+x^2}{2}+\frac{11}{2}+5y^2+e^{2}\\&
\le (7+e^2)(1+|x|^2+|y|^2).\endaligned
\end{equation}

We have shown that condition (\ref{dandiao}) holds for $p=6$ and $K=7+e^2$ for any $x,y.$

Moreover, since $f$ and $g$ are differential on $\mathbb{R}^2$, by mean value theorem, for any $|x|\vee|x'|\vee|y|\vee|y'|\le R,$ let $h=x-\bar{x}$, $k=y-\bar{y}$, then there exists $\theta\in(0,1)$ such that
$$\aligned|f(x,y)-f(x',y')|&=f'_x(\bar{x}+\theta h,\bar{y}+\theta k)h+f'_y(\bar{x}+\theta h,\bar{y}+\theta k)k\\&
\le(2+(1+3R)e^R)(|x-x'|+|y-y'|).\endaligned$$
Similarly,
$$|g(x,y)-g(x',y')|\le Re^R(1+\frac{3}{2}R)(|x-x'|+|y-y'|)$$
for all $R>0$ and $|x|\vee|x'|\vee|y|\vee|y'|\le R.$
Thus (\ref{local}) holds for $L_{R}=3(1+R+R^2)e^R.$

Then for any $0<\varepsilon<1,$ we can define $l(x):=\frac{1}{x^{1-\varepsilon}L^4_x}$ for $x>0$. It is clear that $l$ is a strict decreasing function in the interval $(0,\infty)$. Let $h$ be the inverse function of $l.$ Then $h$ is also a strict decreasing function in the interval $(0,\Delta^*)$ and $h(\Delta)\to\infty$ as $\Delta\to0$.

Now $$L^4_{h(\Delta)}\Delta=L^4_Rl(R)=\frac{1}{R^{1-\varepsilon}},$$
where $R:=h(\Delta)$. Therefore,
$$L^4_{h(\Delta)}\Delta=\frac{1}{h(\Delta)^{1-\varepsilon}}\to0\qquad \textrm{as}\ \Delta\to0.$$

And
$$(L^{2q}_{h(\Delta)}\Delta^{\frac{q}{2}})^{-\frac{1}{p-q}}=(L^{4}_{h(\Delta)}\Delta)^{-\frac{q}{2(p-q)}}
 =h(\Delta)^{\frac{q(1-\varepsilon)}{2(p-q)}}\le h(\Delta)$$
 for $\Delta$ small enough if $\frac{q}{2(p-q)}\le 1$ (i.e. $q\le\frac{2p}{3}=4$). Then by Theorem \ref{conv}, for any $T>0$, $2<q\le4$ and sufficient small $\Delta$, we have
\begin{equation}\label{b1}\mathbb{E}|x(T)-x_\Delta(T)|^q\le C(q,T)L^{2q}_{h(\Delta)}\Delta^\frac{q}{2}=C(q,T)h(\Delta)^{\frac{q(\varepsilon-1)}{2}}\end{equation}
and
\begin{equation}\label{b2}\mathbb{E}|x(T)-\bar{x}_\Delta(T)|^q\le C(q,T)L^{2q}_{h(\Delta)}\Delta^\frac{q}{2}=C(q,T)h(\Delta)^{\frac{q(\varepsilon-1)}{2}}.\end{equation}

H\"{o}lder inequality implies that (\ref{b1}) and (\ref{b1}) holds for any $0<q\le4$.

Since $f$ does not satisfy polynomial growth condition in this case, then the strong convergence result Theorem 3.7 in \cite{ZSL} does not be hold here. However, for the continuous-time MTEM methods (\ref{num1}) and (\ref{num2}), the strong convergence results still holds for the given NSDDE.

\textbf{Example 2}
Consider the scalar NSDDE
\begin{equation}\label{sde2}\aligned d[x(t)-\frac{1}{2}\sin x(t-1)]&=(2x(t)-x^5(t)+\frac{1}{2}\sin
x(t-1))dt+\frac{2x^3(t)x(t-1)}{1+x^2(t-1)}dB_t.\endaligned\end{equation}

Here $f(x,y)=2x-x^5+\frac{1}{2}\sin y$, $g(x,y)=\frac{x^3y}{2(1+y^2)}$ and $D(y)=\frac{1}{2}\sin y.$ It is obvious that $f$ and $g$ are both locally Lipschitz continuous functions with respect to $x$ and $y$.

Moreover, for any $a\in(0,1]$, we have
\begin{equation}
\aligned2\langle
x-aD(\frac{1}{a}y),f(x,y)\rangle+5|g(x,y)|^2&=4x^2-2x^6+x\sin y-2ax\sin \frac{y}{a}+x^5a\sin \frac{y}{a}\\&\quad-\frac{a}{2}\sin y\sin\frac{y}{a}+\frac{5x^6y^2}{4(1+y^2)^2}\\&
\le4x^2-2x^6+|xy|+|2ax\cdot\frac{y}{a}|\\&\quad+\frac{5}{6}x^6+\frac{a^6}{6}\sin^6\frac{y}{a}+\frac{a}{2}|y\cdot\frac{y}{a}|+\frac{5}{16}x^6\\&
\le4x^2+3\frac{x^2+y^2}{2}+\frac{y^2}{6}+\frac{y^2}{2}\\&
\le\frac{11}{2}(1+x^2+y^2).\endaligned
\end{equation}

So (\ref{dandiao}) holds for $p=6$.

Moreover $|g(x,y)|^2\le|x|^3\le 1+|x|^3+|y|^3.$ That is, (\ref{gzeng}) also holds for $r=3$ and $\bar{K}=1$.

On the other hand, we have
$$|f(x,y)-f(x',y')|\le (2+5R^4)(|x-x'|+|y-y'|)$$
and
$$|g(x,y)-g(x',y')|\le(3R^2+R^3)(|x-x'|+|y-y'|)$$
for all $R>0$ and $|x|\vee|x'|\vee|y|\vee|y'|\le R.$
Thus $f$ and $g$ are local Lipschitz continuous with local Lipschitz constant $L_R=5R^4+4.$

For $\varepsilon>0$ small enough, choose
$$h(\Delta)=\sqrt[4]{\frac{\Delta^{-\varepsilon}-4}{5}},\quad \Delta<4^{-\frac{1}{\varepsilon}}.$$
Then we have $h(\Delta)\to \infty$ and $L^4_{h(\Delta)}\Delta=\Delta^{1-\varepsilon}\to0$ as $\Delta\to0.$ That is, (\ref{tiaoj}) holds for such defined $h$.

Choose $q=4$. If we take $\frac{4}{5}<\varepsilon<1,$ then for sufficiently small $\Delta,$
 $$(L^{2q}_{h(\Delta)}\Delta^{\frac{q}{2}})^{-\frac{1}{p-q}}=(L^{4}_{h(\Delta)}\Delta)^{-1}
 =\Delta^{\varepsilon-1}\le \sqrt[4]{\frac{\Delta^{-\varepsilon}-4}{5}}=h(\Delta),$$
 i.e. (\ref{tj}) holds for small $\Delta$. So by Theorem \ref{conv1}, we have
 $$\mathbb{E}\sup_{0\le t\le T}|x(t)-x_\Delta(t)|^4\le C\Delta^{2(1-\varepsilon)},$$
 and
 $$\mathbb{E}\sup_{0\le t\le T}|x(t)-\bar{x}_\Delta(t)|^4\le CL^4_{h(\Delta)}\Delta=C\Delta^{1-\varepsilon}.$$

 \section{Conclusions}

 We have investigated the strong convergence rates of so called two versions of continuous-time MTEM methods (i.e., $x_\Delta(t)$ and $\bar{x}_\Delta(t)$) for nonlinear NSDDE $d[x(t)-D(x(t-\tau))]
=f(x(t),x(t-\tau))dt+g(x(t),x(t-\tau))dB_t$ in this paper. Roughly speaking, $x_\Delta(t)$ and $\bar{x}_\Delta(t)$ strongly converge (in the sense of $q$-th moment) to the exact solution $x(t)$ at fixed time $T$ (with rate $L^{2q}_{h(\Delta)}\Delta^{\frac{q}{2}}$) if local Lipschitz condition and the Khasiminskii-type condition hold. Moreover, if $g$ satisfies polynomial growth condition (\ref{gzeng}), then $x_\Delta(t)$ and $\bar{x}_\Delta(t)$ strongly converge to the exact solution $x(t)$ over a time interval $[0,T]$ (with rates $L^{2q}_{h(\Delta)}\Delta^{\frac{q}{2}}$ and $L^{q}_{h(\Delta)}\Delta^{\frac{q}{2}-1}$, respectively).

\end{document}